\input amstex
\documentstyle{amsppt}
\magnification=\magstep1 \NoRunningHeads
\topmatter

\NoBlackBoxes
\loadbold

\title 
Explicit construction of  orbit~equivalence for rank-one actions
\endtitle

%\comment

\author
Alexandre  I. Danilenko  and Mykyta I. Vieprik
\endauthor

%\endcomment

\abstract
Let $G$ be a discrete countable infinite group. 
Given two  rank-one measure-preserving $G$-actions whose
invariant measures are either both  finite or both infinite, we give a direct explicit
construction, from their  cutting-and-stacking parameters, of a Borel orbit equivalence
on invariant conull subsets.
We also prove a topological counterpart of this assertion, under additional compatibility assumptions, for continuous $(C,F)$-actions of $G$ on non-compact
locally compact Cantor spaces. 
\comment
on standard $\sigma$-finite measure spaces $(X,\mu)$ and $(\widetilde X,\widetilde \mu)$.
If $\mu(X)=\widetilde\mu(X)=1$ or $\mu(X)=\widetilde\mu(X)=\infty$ then  we construct explicitly and algorithmically  a measurable orbit equivalence between  $T$ and $\widetilde T$ in terms of the underlying cutting-and-stacking parameters.
We also prove a topological version of this result is also proved when  $(X,T)$ and $(\widetilde X,\widetilde T)$ are continuous rank-one systems on locally compact non-compact Cantor spaces.

\endcomment
\endabstract

\address
Faculty of Mathematics and Computer Science,
Nicolaus Copernicus University,
Chopin street 12/18, 87-100 Toru{\'n}, Poland
\newline\indent
 B. Verkin Institute for Low Temperature Physics and Engineering
of the  National Academy of Sciences of Ukraine,
47 Nauky Ave.,
 Kharkiv, 61164, Ukraine
\endaddress
\email            alexandre.danilenko\@gmail.com
\endemail

\address
Imperial College London, Department of Mathematics, Huxley Building, 180 Queen’s Gate, London, SW7 2AZ, United Kingdom  
\endaddress
\email
nikita.veprik\@gmail.com
\endemail

\endtopmatter
\document

\head 0. Introduction
\endhead

Orbit equivalence ignores the time parametrization of an action and preserves only the partition into orbits. 
Therefore, classifying actions up to orbit equivalence is much simpler than classifying them up to isomorphism.
A fundamental result of Dye  asserts that any  
two ergodic  measure-preserving automorphisms of a nonatomic standard probability space are orbit equivalent \cite{Dy};
an infinite measure analog of this theorem is due to Krieger  \cite{Kr}.
%Similar result holds  in the infinite $\sigma$-finite measure spaces \cite{Kr}.
For extension of these results  to nonsingular actions of discrete countable amenable groups we refer to \cite{OrWe} and  \cite{CoFeWe}.
%In the case of continuous measure preserving $\Bbb Z$-actions on Polish spaces, Dye's theorem admits a refinement: the measurable orbit equivalence can be chosen to be almost continuous,
%that is, implemented by a homeomorphism between dense $G_\delta$ subsets of full measure
%(see \cite{HaKe}, \cite{dJSa}, \cite{DadJ} and references therein).
However, the general constructions of  orbit equivalences in the literature are implicit:
they do not provide a recipe for turning the parameters determining the actions into
a concrete orbit equivalence map admitting
a finite-stage description.
 The goal of this paper is to propose such a parametric construction for the {\it rank-one} actions.

%Our main objective in this paper is to  develop an explicit and algorithmic construction of  orbit equivalences for a broad class of actions, namely the rank-one actions. 
Let $G$ be an infinite countable discrete group.
There exist several  definitions (and hence models) of  rank-one $G$-actions: see, for example,  \cite{Fe}  for the case  $G=\Bbb Z$.
Our approach is based on a model   realized via an inductive cutting-and-stacking process  governed by  a sequence 
$\Cal T=(C_n,F_{n-1})_{n=1}^\infty$
of finite subsets of $G$ satisfying the usual cutting-and-stacking conditions.
We call $\Cal T$  the $(C,F)$-parameters.
 %The resulting systems will be referred to as $(C,F)$-actions (\cite{Da1}, 
 %\cite{Da3}).
 The associated system, called the $(C,F)$-action, is a concrete rank-one model on a totally disconnected Polish space; conversely, every rank-one action admits a measure-theoretic  $(C,F)$-model \cite{Da3}, \cite{DaVi1}.
 The advantage of this realization is that telescopings, reductions and the orbit maps constructed below are all described by finite operations on the sets $(C_n)_{n=1}^\infty$ and $(F_n)_{n=0}^\infty$.
 
 Let $T=(T_g)_{g\in G}$ and $\widetilde T=(\widetilde T_g)_{g\in G}$ be two measure-preserving $G$-actions on standard measure spaces. 
An orbit equivalence between them  is a nonsingular Borel isomorphism between invariant conull Borel subsets which maps each $T$-orbit onto a $\widetilde T$-orbit. 
If $T$ and $\widetilde T$ are ergodic then the orbit equivalence is measure scaling.

 \comment
 
 Thus, the $(C,F)$-actions are nothing but the rank-one actions realized in a specific way determined by $\Cal T$.
 They possess an additional topological structure:
 they appear as
  continuous, minimal, free  amenable actions\footnote{$G$ itself is not assumed amenable, in general.}
on totally disconnected Polish spaces equipped with either a finite or an infinite invariant measure. 
An arbitrary  $G$ admits a  rank-one action with an infinite $\sigma$-finite invariant measure.
If this measure is finite then $G$ is amenable.
\footnote{It remains an open question whether every amenable group admits a rank-one action with finite invariant measure.}

\endcomment

Our first main result gives a direct measurable construction\footnote{For  rigorous definitions of the $(C,F)$-space $X$,  the Haar measure $\mu$  and the $(C,F)$-action $T$  associated with $\Cal T$ see \S1.1 and \S1.2.}.

\proclaim{Theorem A} Let $G$ be a discrete countable infinite  group and let $T=(T_g)_{g\in G}$ and $\widetilde T=(\widetilde T_g)_{g\in G}$  be
 two  measure-preserving $(C,F)$-actions on the $(C,F)$-spaces $(X,\mu)$ and $(\widetilde X,\widetilde\mu)$ associated with $(C,F)$-parameters $\Cal T$ and $\widetilde{\Cal T}$ respectively
 such that either 
 $$
\max( \mu(X),\widetilde\mu(\widetilde X))<\infty\quad\text{or}\quad \min(\mu(X),\widetilde\mu(\widetilde X))=\infty.
 $$
 Then there exist  conull invariant Borel subsets $X^\circ\subset X$ and $\widetilde X^\circ\subset \widetilde X$, a Borel bijection $\phi:X^\circ\to \widetilde X^\circ$ and $c>0$ such that  
 $$
 \widetilde\mu\circ\phi=c\cdot \mu\quad\text{and}\quad
 \phi(\{T_gx\mid g\in G\})=\{\widetilde T_g\phi(x)\mid g\in G\}\quad\text{for each $x\in X^\circ$.}
 $$
 If $\max(\mu(X),\widetilde\mu(\widetilde X))<\infty$ then $c=\widetilde\mu(\widetilde X)/\mu(X)$.
The objects   $X^\circ$, $\widetilde X^\circ$ and  $\phi$ are constructed  explicitly   in terms of $\Cal T$ and $\widetilde{\Cal T}$.
\endproclaim

The novelty of Theorem A is in the final claim.
We use   ``explicitly   in terms of $\Cal T$ and $\widetilde{\Cal T}$\,'' in a finite-stage sense: after enumerations of the finite parameter sets have been fixed, the subsets, injections and bijections entering each stage of the inductive construction 
are specified by finite searches. 
The construction has two ingredients. 
First, suitable telescopings and  reductions of the two parameter sequences are chosen so that their new cut sets $C_n'$ and $\widetilde C_n'$ admit alternating factorizations
$$
\# C_n'=a_nb_n,
\qquad
\# \widetilde C_n'=b_na_{n+1}.
$$
Second, an explicit coordinate interchange, together with finite injections between the tower shapes, produces an orbit equivalence of the reduced systems.
Thus, the argument makes the orbit equivalence visible at every finite stage.
%No recursion-theoretic assertion about arbitrary real input data is intended.

%We illustrate Theorem A by an  explicit construction of an  orbit equivalence between the 2-adic and 3-adic odometers in Example~3.4 (cf. \cite{HaKe}).

We also investigate a topological counterpart of this problem. 
Under suitable conditions on $\Cal T$ (see \cite{Da3, Proposition~2.2}), 
the space of a $(C,F)$-action is not only Polish but also a non-compact, locally compact Cantor space equipped with a Radon invariant measure, and the action is  continuous, minimal and Radon uniquely ergodic.
For $G=\Bbb Z$, the topological orbit equivalence
in this setting is governed by ordered cohomological and dimension-group invariants
 \cite{Da2}, \cite{Ma}.
 The special form of a $(C,F)$-system leads to a simple parameter invariant. 
Given a $(C,F)$-sequence $\Cal T=(C_n,F_{n-1})_{n=1}^\infty$, we put
$$
p_0(\Cal T):=1,\qquad p_n(\Cal T):=\#C_1\cdots\#C_n\quad(n\ge1),
$$
and
$$
\Cal D(\Cal T):=\bigcup_{n\ge0}\frac1{p_n(\Cal T)}\Bbb Z.
$$
The following statement  is the second main result of the paper.
%concrete compatibility condition yields an explicit orbit equivalence.

\proclaim{Theorem B} 
Let $T=(T_g)_{g\in G}$ and $\widetilde T=(\widetilde T_g)_{g\in G}$ be
 two continuous $(C,F)$-actions of $G$ on non-compact locally compact Cantor spaces $X$ and $\widetilde X$, associated with $(C,F)$-sequences $\Cal T$ and $\widetilde{\Cal T}$ respectively.
 Let $\mu$ and $\widetilde\mu$ be the corresponding invariant 
 Haar measures.  
 Then 
 there 
 exists a homeomorphism $\phi:X\to \widetilde X$ such that
 $$
 \phi(\{T_gx\mid g\in G\})=\{\widetilde T_g\phi(x)\mid g\in G\}
 \quad\text{for each $x\in X$}
 $$
%and $\mu\circ\phi^{-1}=c\cdot\widetilde\mu$
if and only if 
 there is $c>0$
with
 $$
 \Cal D(\Cal T)= c\cdot\Cal D(\widetilde{\Cal T})\qquad\text{and}\qquad
 \mu(X)=c\widetilde\mu(\widetilde X).
 $$
Moreover, in the ``if'' case,
$\phi$ is constructed explicitly and algorithmically 
 in terms of $\Cal T$, $\widetilde{\Cal T}$ and  $c$.
 Furthermore, $\mu\circ\phi^{-1}=c\widetilde\mu$.
 %the underlying $(C,F)$-parameters.
\endproclaim

%The extra  hypothesis 
%$\Cal D(\Cal T)=c\cdot\Cal D(\widetilde{\Cal T})$
%in Theorem~B cannot be omitted.  
%Proposition~2.5 below gives a compact-open measure-value obstruction, and 
Example~2.5 gives two continuous rank-one $\Bbb Z$-actions on non-compact locally compact Cantor spaces which are not topologically orbit equivalent.

Our approach to Theorems~A and B builds on modifications of techniques developed in
 \cite{DaVi2} and \cite{Da4} for the classification of $(C,F)$-actions up to measurable and topological isomorphism, respectively.
 
 Section~1 recalls the $(C,F)$-construction, telescopings, reductions and the rank-one model theorem. 
 Section 2 proves the coordinate-interchange theorem (Theorem~2.1), 
 derives Theorem~B and records the compact-open measure-value obstruction (Example~2.5). 
 Section 3 proves Theorem~A.

\head 1. $(C,F)$-construction
\endhead

\subhead 1.1. $(C,F)$-spaces, tail equivalence relations and  return time cocycles 
\endsubhead
We remind the $(C,F)$-construction as it appeared in \cite{Da3}.
Let $G$ be  a discrete countable group.
Let $\Cal T=(C_n,F_{n-1})_{n=1}^\infty$ be  a sequence of (pairs of) finite subsets  of $G$ such that
$\#F_{0} = 1$ and for 
 each $n>0$,
$$
\aligned
	  &  
	  \#C_{n} > 1, \\ 
      &F_{n} C_{n+1}\subset F_{n+1},\\ 
      &F_{n} c\cap F_{n} c' = \emptyset\text{ if $c, c'\in C_{n+1}$ and $c \neq c'$. }
\endaligned
\tag1-1
$$
We let $X_n := F_{n} \times C_{n+1} \times C_{n+2} \times\ldots$ and endow this set with the infinite product topology. 
Then $X_n$ is a compact Cantor space. The mapping 
$$
	X_n \ni (f_n,c_{n+1},c_{n+2}\ldots) \mapsto (f_n c_{n+1},c_{n+2},\ldots) \in X_{n+1}\tag1-2
$$
is a continuous embedding of $X_n$ into $X_{n+1}$. 
Therefore, the topological  inductive limit $X$ of the sequence $(X_n)_{n\geq 0}$  is  well defined.
Of course, $X$ is a locally compact Cantor space. 
We call $X$ {\it the $(C,F)$-space associated with $\Cal T$}.
It is convenient to think of  $X$ as the union $\bigcup_{n=1}^\infty X_n$ of the increasing
sequence $X_0\subset X_1\subset\cdots$ of compact open subsets, where the corresponding embeddings are given by \thetag{1-2}.
For a subset $A\subset F_n$, we let
$$
	[A]_n := \{x=(f_n,c_{n+1},\ldots)\in X_n, f_n\in A\}
$$
and call this set an $n$-{\it cylinder} in $X$.
It is open and compact in $X$. 
Every open subset of $X$ is a union of cylinders.
For brevity, we will write $[f]_n$ for $[\{f\}]_n$ for an  element $f\in F_n$.

Two points $x=(f_n,c_{n+1},\ldots)$ and $x'=(f_n',c_{n+1}',\ldots)$ of $X_n$ are called {\it tail equivalent}
if there is $N>n$ such that $c_l=c_{l}'$ for  each $l>N$.
We thus obtain {\it the tail equivalence relation} $\Cal R_n$ on $X_n$.
The {\it  tail equivalence relation}  $\Cal R$ on $X$ is defined as follows: for each $n\ge 0$, the restriction of $\Cal R$ to $X_n$ is $\Cal R_n$.
It is well defined and Borel.
Of course, $\Cal R$ is {\it minimal} on $X$, i.e. each $\Cal R$-class is dense in $X$.

We note that $\Cal R$ is {\it Radon uniquely ergodic}, i.e. there is a unique  $\Cal R$-invariant Radon measure $\mu$ on $X$ such that $\mu(X_0)=1$.
We call it  the {\it  Haar measure  for $\Cal R$}.
The Haar measure is $\sigma$-finite.
Let $\kappa_n$ denote the equidistribution on $C_n$.
We define a measure $\nu_n$ on $F_n$ by setting 
$$
\nu_0(F_0):=1\quad\text{and \,
 $\nu_n(\{f\}):=
 \frac1{\prod_{k=1}^n \#C_k}$ for each  $f\in F_n$ and $n> 0$.}
 $$
 Then
 $$
 \mu\restriction X_n=\nu_n\otimes\bigotimes_{k>n}\kappa_k\quad\text{for each $n>0$}.
 $$
We note that $\mu$ is finite if and only if 
$$
\prod_{n=1}^\infty\frac{\# F_{n+1}}{\# F_n\# C_{n+1}}<\infty.\tag1-3
$$

Define a Borel mapping $\alpha:\Cal R\to G$ by setting
$$
\alpha(x,\widetilde x):=\lim_{m\to\infty}f_nc_{n+1}\cdots c_m\widetilde c_m^{-1}\cdots\widetilde c_{n+1}^{-1}\widetilde f_n^{-1}
$$
if $x=(f_n,c_{n+1},c_{n+2},\dots)\in X_n$ and 
$\widetilde  x=(\widetilde f_n,\widetilde c_{n+1},\widetilde c_{n+2},\dots)\in X_n$
for some $n\ge 0$.
It is straightforward to verify that  $\alpha$ is well defined and it satisfies the cocycle identity
$$
\alpha(x,\widetilde x)\alpha(\widetilde x,\widehat x)=\alpha(x,\widehat x)
$$ 
for all $x,\widetilde x,\widehat x\in X$
such that $(x,\widetilde x)\in\Cal R$ and $(\widetilde x,\widehat x)\in\Cal R$.
We call $\alpha$ {\it the return time cocycle} of $\Cal R$.

\subhead 1.2. $(C,F)$-actions: topological and measure-theoretic
\endsubhead
Given $g\in G$,  let 
$$
X_n^g:=\{(f_n, c_{n+1}, c_{n+2},\dots)\in X_n\mid gf_n\in F_n\}.
$$
Then $X_n^g$ is a compact open subset of $X_n$ and $X_n^g\subset X^g_{n+1}$.
  Hence the union
  $X^g:=\bigcup_{n\ge 0}X_n^g$ is an open subset of $X$.
Let 
$$
X^G:=\bigcap_{g\in G}X^g.
$$
 Then $X^G$ is a $G_\delta$-subset of $X$.
  Hence $X^G$ is Polish and totally disconnected  in the induced topology. 
  Given $g\in G$ and  $x\in X_G$, there is $n >0$ such that
  $x= (f_n, c_{n+1},\dots )\in X_n$ and $gf_n\in F_n$. 
  We now let 
  $$
  T_gx:= (gf_n, c_{n+1}, \dots)\in X_n\subset X.
  $$
  It is straightforward  to verify that
  \roster
  \item"---" $T_gx\in X^G$,
  \item"---" the mapping $T_g:X^G\ni x\mapsto T_gx\in X^G$ is a well defined homeomorphism of $X^G$, 
  \item"---" $T_gT_{g'}=T_{gg'}$ for all $g, g'\in G$ and 
  \item"---"  $\alpha(T_gx,x)=g$ for all $g\in G$ and $x\in X^G$,
  where $\alpha$ is the return time cocycle of $\Cal R$.
\endroster
Hence, $T:=(T_g)_{g\in G}$
is a continuous, well defined $G$-action on $X^G$.

\definition{Definition 1.2 \cite{Da3}} The action $T$ is called {\it the topological $(C,F)$-action of $G$ associated with $\Cal T$}.
\enddefinition

The topological $(C,F)$-action is free.
The subset $X^G$ is $\Cal R$-invariant.
The $T$-orbit equivalence relation 
coincides with the restriction of $\Cal R$ to $X^G$.
It was shown in \cite{Da3, Proposition~2.2} that
$X^G=X$ if and only if for each $g\in G$ and $n >0$, there is
$m > n$ such that
$$
gF_nC_{n+1}C_{n+2}\cdots C_m\subset F_m.
\tag1-4
$$
Thus, if \thetag{1-4} holds then $T$ is a minimal continuous $G$-action on a locally compact Cantor space $X$.
This action is {\it Radon uniquely ergodic}, i.e. there exists a unique  $T$-invariant Radon measure $\xi$ on $X$ such that $\xi(X_0)=1$.
Of course, $\xi$ is the Haar measure for  $\Cal R$.
If  $T$ is Radon uniquely ergodic and \thetag{1-3} holds then $T$ is uniquely ergodic in the classical sense, i.e. there exists a unique invariant probability Borel measure on $X$.

From now on, $T$ is a topological $(C,F)$-action of $G$ on $X^G$ and
 $\mu$ is the Haar measure for $\Cal R$. 
 Since $X^G$ is $\Cal R$-invariant, we obtain that either $\mu(X^G)=0$ or $\mu(X\setminus X^G)=0$.
In the latter case $T$ is  conservative and ergodic.
The following two results were obtained  in \cite{Da3, Proposition~2.4}:

\proclaim{Fact A}
$\mu(X\setminus X^G)=0$ if and only if for each $g\in G$ and every $n \ge0$,
$$
\lim_{m\to\infty}\nu_m\big((gF_nC_{n+1}C_{n+2}\cdots C_m)\cap F_m\big)=\nu_n(F_n).\tag1-5
$$
\endproclaim

\proclaim{Fact B}
If  $\mu(X\setminus X^G)=0$ and $\mu(X)<\infty$  then
$G$ is amenable and 
$(F_n)_{n=1}^\infty$ is a left F{\o}lner sequence in $G$.
\endproclaim

\comment
We also note that if $G$ admits a finite measure-preserving $(C,F)$-action then
the F{\o}lner sequence $(F_n)_{n=1}^\infty$ possesses the following ``near tiling'' property:  for each pair  of integers $m>n>0$, there is a finite subset $D_{n,m}$ such that $F_nD_{n,m}\subset F_m$, $F_nd\cap F_nd'=\emptyset$ for all $d\ne d'\in D_{n,m}$ and for each $\epsilon>0$, if $n$ is large enough then
$\#(F_nD_{n,m})/\#F_m>1-\epsilon$ for every $m>n$.
We do not know whether each countable amenable group has a F{\o}lner sequence with the near tiling property.
\endcomment

\definition{Definition 1.3} If $\mu(X\setminus X^G)=0$ then the dynamical system $(X,\mu, T)$ (or simply $T$) is called
{\it the  measure-preserving $(C,F)$-action associated with  $\Cal T$.}
\enddefinition

\subhead  1.3. Telescopings 
\endsubhead
Let a sequence $\Cal T=(C_n,F_{n-1})_{n=1}^\infty$ satisfy   \thetag{1-1}. 
Given a strictly increasing  infinite sequence of integers $\boldsymbol l=(l_n)_{n=0}^\infty$ such that
$l_0=0$, we let 
$$
\widetilde F_n:=F_{l_n},\quad \widetilde C_{n+1}:=C_{l_n+1}\cdots C_{l_{n+1}}
$$
for each $n\ge 0$.
The sequence $\widetilde{\Cal T}:=(\widetilde C_n,\widetilde F_{n-1})_{n=1}^\infty$ is called the {\it $\boldsymbol l$-tele\-scoping} of $\Cal T$ \cite{Da3}.
It is easy to check that $\widetilde{\Cal T}$   satisfies  \thetag{1-1}.  
Let $X$ and $\widetilde X$ denote the $(C,F)$-spaces associated with $\Cal T$
and $\widetilde{\Cal T}$ respectively.
Denote by $\Cal R$ and $\widetilde{\Cal R}$ the tail equivalence relations on $X$ and $\widetilde X$
respectively.
There is a canonical mapping $\iota_{\boldsymbol l}$ of  $X$ onto $\widetilde X$ associated with $\boldsymbol l$.
If $x\in X$ then we select the smallest  $n\ge 0$ such that $x=(f_{l_n},c_{l_n+1},c_{l_n+2},\dots)\in X_{l_n}$ and put
$$
\iota_{\boldsymbol l}(x):=(f_{l_n},c_{l_n+1}\cdots c_{l_{n+1}}, c_{l_{n+1}+1}\cdots c_{l_{n+2}},\dots)\in \widetilde X_n\subset\widetilde X,
$$
where $\widetilde X_n=\widetilde F_n\times\widetilde C_{n+1}\times \widetilde C_{n+2}\times\cdots$.
It is routine to verify that 
\roster
\item"---" $\iota_{\boldsymbol l}$ is a homeomorphism of $X$
onto $\widetilde X$,
\item"---" $\iota_{\boldsymbol l}$ maps bijectively each $\Cal R$-class in $X$ onto
an $\widetilde{ \Cal R}$-class in $\widetilde X$,
\item"---" 
$\iota_{\boldsymbol l}$ transfers the Haar measure for $\Cal R$ to the Haar measure
for $\widetilde {\Cal R}$.
\endroster
Moreover,
$$
\widetilde\alpha(\iota_{\boldsymbol l}x,\iota_{\boldsymbol l}x')=\alpha(x,x')\quad\text{for each $(x,x')\in\Cal R$},\tag{1-6}
$$
where $\alpha$ and $\widetilde\alpha$ are the return time cocycles of $\Cal R$ and $\widetilde{\Cal R}$ respectively.
We call $\iota_{\boldsymbol l}$ {\it the $\boldsymbol l$-telescoping mapping.}

If $\Cal T$ satisfies  \thetag{1-4} or \thetag{1-5} then 
$\widetilde{\Cal T}$ also satisfies~\thetag{1-4} or \thetag{1-5} respectively.
Hence,  the $(C,F)$-actions $T$ and $\widetilde T$ associated with $\Cal T$ and $\widetilde{\Cal T}$ respectively
are well defined.
It follows from  \thetag{1-6} that $T$ and $\widetilde T$ are conjugate via $\iota_{\boldsymbol l}$, i.e.
$\iota_{\boldsymbol l}T_g\iota_{\boldsymbol l}^{-1}=\widetilde T_g$ for each $g\in G$.

\subhead 1.4. Reductions
\endsubhead
Let a sequence $\Cal T=(C_n,F_{n-1})_{n=1}^\infty$ satisfy   \thetag{1-1}.
Let $\boldsymbol A:=(A_n)_{n=1}^\infty$ be  a sequence of nonempty  subsets $A_n\subset C_{n}$ such that $\# A_n\ge 2$ for each $n\in\Bbb N$ and 
$$
\sum_{n=1}^\infty(1- \kappa_n(A_n))<\infty.
$$ 
We will assume that $A_n$ is a proper subset of $C_n$ for infinitely many $n$.
Denote by $ \kappa_n^*$
the equidistribution 
on $A_n$ for each $n\in\Bbb N$.
Let $\Cal T^*:=(A_n,F_{n-1})_{n=1}^\infty$.
The sequence $\Cal T^*$  is called   {\it the $\boldsymbol A$-reduction} of $\Cal T$
\cite{Da3}.
It is easy to check that $\Cal T^*$   satisfies  \thetag{1-1}.
Let $X$ and $X^*$ be the $(C,F)$-spaces associated with $\Cal T$ and $\Cal T^*$.
Denote by 
$\Cal R$ and $\Cal R^*$ the tail equivalence relations on $X$ and $ X^*$
respectively.
Let $\mu$ and $\mu^*$ denote the Haar measures on $X$ and $X^*$
respectively.
We note that for each $n\ge0$, the identity mapping embeds the set 
$$
X_n^*:=F_n\times A_{n+1}\times A_{n+2}\times\cdots
$$
into  $X_n:=F_n\times C_{n+1}\times C_{n+2}\times\cdots$.
Hence, we can think of  $X_n^*$ 
as a nowhere dense closed subset of $X_n$.
It follows that $X^*=\bigcup_{n\ge 0}X_n^*$ embeds naturally  into $X$ as an $F_\sigma$-subset of the first Baire category.
Of course, $X^*$ is $\Cal R$-invariant and, hence, dense in $X$.
The restriction of $\Cal R$ to $X^*$ is $\Cal R^*$.
We note that
$$
\mu(X^*)\ge \mu(X_n^*)\ge \prod_{j>n}{\kappa_j(A_j)}>0.
$$
Since $X^*$ is $\Cal R$-invariant and $\mu$ is $\Cal R$-ergodic, it follows that $\mu(X\setminus X^*)=0$.
Thus, $X^*$ is of full Haar  measure in $X$.
There is a  canonical measure scaling Borel isomorphism $\rho_{\boldsymbol A}$ of  $(X,\mu)$ onto $(X^*,\mu^*)$:
$$
\rho_{\boldsymbol A}x:=x\quad\text{if $x\in X^*$}.
$$
The reader should not confuse $x$ from the left-hand side ($x$ is a point of the $(C,F)$-space $X$) with $x$ from the right-hand side ($x$ is a point of the $(C,F)$-space $X^*$ which carries a different topology).
Thus, $\rho_{\boldsymbol A}$ is defined only on $X^*$, which is a
$\mu$-conull  $F_\sigma$- subset of  $X$.
It is straightforward to verify that
\roster
\item"---" the inverse mapping $\rho_{\boldsymbol A}^{-1}:X^*\to X$ 
is continuous, 
\item"---"
$\rho_{\boldsymbol A}$ maps bijectively each $\Cal R\cap(X^*\times X^*)$-class in $X^*\subset X$ onto
an $\Cal R^*$-class in $X^*$,
\item"---"
$
\frac{d(\mu\circ \rho_{\boldsymbol A}^{-1})}{d\mu^*}=\prod_{m>0}{\kappa_m(A_m)}
$
almost everywhere
and
\item"---"  
$\alpha^*(\rho_{\boldsymbol A}x,\rho_{\boldsymbol A}x')=\alpha(x,x')$ for each $(x,x')\in\Cal R\cap(X^*\times X^*)$,
\endroster
where $\alpha$ and $\alpha^*$ are the return time cocycles of $\Cal R$ and $\Cal R^*$ respectively.
We call $\rho_{\boldsymbol A}$ {\it the $\boldsymbol A$-reduction mapping.}

If $\Cal T$ satisfies  \thetag{1-4} or \thetag{1-5} then 
$\Cal T^*$ also satisfies~\thetag{1-4} or \thetag{1-5} respectively.
Hence,  the $(C,F)$-actions $T$ and $ T^*$ associated with $\Cal T$ and $\Cal T^*$ respectively
are well defined.
Moreover,  $T$ and $T^*$ are conjugate via $\rho_{\boldsymbol A}$, i.e.
$\rho_{\boldsymbol A}T_g\rho_{\boldsymbol A}^{-1}= T_g^*$ a.e. for each $g\in G$.

\proclaim{Fact C \cite{Da3}} Let $\Cal T$ be a $(C,F)$-sequence satisfying \thetag{1-1}, \thetag{1-3} and \thetag{1-5}.
Then there is a telescoping $\widetilde {\Cal T}$ of $\Cal T$ 
and a reduction  $\widetilde {\Cal T}^*$ of  \, $\widetilde {\Cal T}$ such that 
$\widetilde {\Cal T}^*$ satisfies~\thetag{1-1}, \thetag{1-3} and \thetag{1-4}.
\endproclaim

\subhead 1.5. Actions of rank one
\endsubhead
 Let $T= (T_g)_{g\in G}$ 
 be a free measure-preserving action of $G$ on a standard $\sigma$-finite non-atomic measure space $(X,\goth B,\mu)$.
  By a {\it Rokhlin tower for} $T$ we mean a pair $(B,F)$, where $B\in\goth B$ 
  with $0<\mu(B)<\infty$ and $F$ is a finite subset of $G$  such that $1_G\in F$ and
  the subsets $T_fB$, $f\in F$, are mutually disjoint.
  Given a Rokhlin tower $(B,F)$, we let $X_{B,F}:=\bigsqcup_{f\in F}T_fB$.
  By $\xi_{B,F}$ we mean the finite partition of $X_{B,F}$ into the subsets $T_fB$, $f\in F$.
  If $x\in T_fB$ then we set $O_{B,F}(x) :=\{T_gx\mid g\in Ff^{-1}\}$.
  Let $\{1_G\}=F_0\subset F_1\subset F_2\subset\dots$ be an increasing sequence of finite subsets in $G$.

 We say that $T$ is {\it of rank one along} $(F_n)_{n=0}^\infty$ if there is a decreasing sequence
   $B_0\supset B_1\supset B_2\supset\cdots $ of subsets of positive measure in $X$ such that
   $(B_n,F_n)$ is a Rokhlin tower for $T$ for each $n\in\Bbb N$ and
   \roster
   \item"---" $X_{B_1, F_1}\subset X_{B_2, F_2}\subset\cdots$, $\bigcup_nX_{B_n, F_n}=X$,
   \item"---" $\xi_{B_1,F_1}\prec\xi_{B_2,F_2}\prec\cdots$, $\bigvee_n\xi_{B_n,F_n}$ is the partition of $X $ into singletons and 
   \item"---"
   $\{T_gx\mid g\in G\}=\bigcup_{n}O_{B_n,F_n}(x)$
   \endroster
at a.e. $x\in X$.

\proclaim{Fact D \cite{DaVi1}} If $T$ is a $(C,F)$-action then $T$ is of rank one. 
Each rank-one action is measure-theoretically isomorphic to a $(C,F)$-action. 
\endproclaim

If $G=\Bbb Z$ then the classical finite measure-preserving transformations of rank-one correspond to  the rank one along a subsequence of the sequence $(\{0,1,\dots,n\})_{n=0}^\infty$.

It follows from Facts~C and D that  each finite measure-preserving rank-one action of an arbitrary $G$
 is measure-theoretically isomorphic 
to a minimal uniquely ergodic continuous $(C,F)$-action on a locally compact Cantor space.

\head 2. Topological orbit equivalence and proof of Theorem B
\endhead

  Let $G$ be a discrete countable infinite group. 
 
\proclaim{Theorem 2.1} Let $(X,\mu)$ and $(\widetilde X,\widetilde\mu)$ be the $(C,F)$-measure spaces associated with 
$(C,F)$-sequences
$\Cal T:=(C_n,F_{n-1})_{n=1}^\infty$ and $\widetilde{\Cal T}:=(\widetilde C_n,\widetilde F_{n-1})_{n=1}^\infty$ respectively.
Let $\Cal R$ and $\widetilde{\Cal R}$ be the corresponding tail equivalence relations on 
$X$ and $\widetilde X$ respectively.
Suppose that there exist two sequences $(a_n)_{n=1}^\infty$ and $(b_n)_{n=1}^\infty$ of natural numbers
such that  $a_1=1$ and for each $n\in\Bbb N$,
\roster
\item"(i)" 
 $\#C_n=a_nb_n$ and $\#\widetilde C_n=b_na_{n+1}$ and
\item"(ii)"  $\widetilde\mu(\widetilde X_{n-1})\le \mu( X_{n})\le \widetilde\mu(\widetilde X_{n})$.
\endroster
Then there exists a homeomorphism $\phi:X\to\widetilde X$ such that $(\phi\times\phi)(\Cal R)=\widetilde{\Cal R}$ and $\mu\circ\phi^{-1}=\widetilde \mu$.
The homeomorphism is constructed  
in an explicit algorithmic way  from $\Cal T$ and $\widetilde{\Cal T}$.
\endproclaim

\demo{Proof} For each $n\in\Bbb N$,  we select sets $A_n$ and $B_n$ of cardinality $a_n$ and $b_n$ respectively and fix 
 a bijection between $C_n$ and  $A_n\times B_n$ and a bijection between $\widetilde C_n$ and  $B_n\times A_{n+1}$.
 Such bijections exist due to (i).
 For simplicity, we will suppress the notation of these bijections and simply write an element $c\in C_n$ as a pair $c=(a,b)$ with $a\in A_n$ and $b\in B_n$, and an element
 $\widetilde c\in \widetilde C_n$ as a pair $\widetilde c=(d,e)$ with $d\in B_n$ and $e\in A_{n+1}$.
 
  It follows from (ii) and (i) that
 $$
 1\le \frac{\widetilde\mu(\widetilde X_{n})}{\mu(X_{n})}=\frac{\#\widetilde F_n}{\#\widetilde C_1\cdots\#\widetilde C_n}\cdot
 \frac{\#C_1\cdots\#C_n}{\#F_n}=\frac{\#\widetilde F_n}{\#F_n\cdot a_{n+1}}.
 $$
 Thus,  ${\#\widetilde F_n}\ge\#F_n\cdot \# A_{n+1}$.
 In a similar way, we deduce from (ii) and (i) that
 $$
 1\ge\frac{\widetilde\mu(\widetilde X_{n-1})}{\mu(X_{n})}=\frac{\#\widetilde F_{n-1}}{\#\widetilde C_1\cdots\#\widetilde C_{n-1}}\cdot
 \frac{\#C_1\cdots\#C_n}{\#F_n}=\frac{\#\widetilde F_{n-1}\cdot b_n}{\#F_n},
 $$
 i.e. $\# F_n\ge \#\widetilde F_{n-1}\cdot \#B_n$.
 Hence, 
 $$
 \# F_n\cdot\# A_{n+1}\ge \#\widetilde F_{n-1}\cdot \#B_n\# A_{n+1}=\#\widetilde F_{n-1}\cdot\#\widetilde C_n.
 $$
 Thus, we obtain a double inequality
 $$
\#\widetilde F_n \ge \# F_n\cdot\# A_{n+1}\ge
 \#\widetilde F_{n-1}\cdot\#\widetilde C_n.\tag2-1
 $$
 We now define  inductively  a sequence of  mappings $\phi_n:F_n\times A_{n+1}\to\widetilde F_n$, $n=0,1,\dots$, such that
 
 \roster
 \item"$(\circ)$" $\phi_n$ is one-to-one,
 \item"$(\bullet)$"
 $\phi_n (F_n\times A_{n+1})\supset  \widetilde F_{n-1}\widetilde C_n$,
 \item"$(\ast)$"
 for each $f\in F_{n-1}$, $c=(a,b)\in A_n\times B_n=C_n$ and $a_{n+1}\in A_{n+1}$, we have that
 $$
 \phi_n (fc,a_{n+1})=\phi_{n-1}(f, a)\widetilde c,\tag2-2
 $$ 
 where $\widetilde c:=(b, a_{n+1})\in B_n\times A_{n+1}=\widetilde C_n$.
 \endroster
Of course, $\phi_0$ is defined in an obvious (and unique) way.
Suppose that  $\phi_1,\dots,\phi_{n-1}$ have been already defined for some $n$.
Our purpose is to define $\phi_n$.
This will be done in two steps.
First, we define $\phi_n$ on the subset $(F_{n-1}C_n)\times A_{n+1}$ via~\thetag{2-2}.
This mapping is one-to-one and its image equals
$\phi_{n-1}(F_{n-1}\times A_{n})\widetilde C_{n} $ which is a subset of 
$\widetilde F_{n-1}\widetilde C_{n}$.
 The ``left'' inequality in~\thetag{2-1} implies that $\phi_n$ extends as a one-to-one mapping
from the entire set
$F_n\times A_{n+1}$  to $\widetilde F_{n}$.
Moreover, in view of the ``right'' inequality in~\thetag{2-1}, we can choose the extension so that 
$(\bullet)$ holds.
Thus, $\phi_n$ is as desired.

We now define  a mapping $\phi:X\to\widetilde X$.
Take $x\in X$.
Then there is $n\ge 0$ such that $x\in X_n$.
Hence, we can write $x$ as a sequence $x=(f_n,c_{n+1},c_{n+2},\dots)$.
Then for each $m>n$, there are $a_m\in A_m$ and $b_m\in B_m$ such that $c_m=(a_m,b_m)$.
Thus,
$$
x=(f_n, (a_{n+1}, b_{n+1}), (a_{n+2}, b_{n+2}),\dots).
$$
We now set
$$
\phi(x):=(\phi_n(f_n,a_{n+1}), (b_{n+1},a_{n+2}), (b_{n+2},a_{n+3}),\dots)\in \widetilde X_n\subset\widetilde X.
$$

{\sl Claim 1.}   $\phi$ is well defined.

\noindent Indeed, we can write $x$ in a different way: $x=(f_n\cdot (a_{n+1}, b_{n+1}),(a_{n+2}, b_{n+2}),\dots )\in X_{n+1}$.
Then, in view of $(\ast)$,
$$
\align
\phi(x)&=
(\phi_{n+1}(f_n\cdot (a_{n+1}, b_{n+1}), a_{n+2}), (b_{n+2},a_{n+3}),\dots)\\
&=(\phi_n(f_n,a_{n+1})\cdot
(b_{n+1},a_{n+2}), (b_{n+2},a_{n+3}),\dots)\\
&=(\phi_n(f_n,a_{n+1}),
(b_{n+1},a_{n+2}), (b_{n+2},a_{n+3}),\dots),
\endalign
$$
which proves  the claim.  

{\sl Claim 2.}   $\phi$ is one-to-one.

\noindent This claim follows from the definition of $\phi$ and $(\circ)$ in a routine way.

{\sl Claim 3.}   $\phi$ is onto.

\noindent It follows  from $(\bullet)$ that for each $n\in\Bbb N$,
 $$
 \widetilde X_{n-1}= [\widetilde F_{n-1}\widetilde C_n]_n\subset [\phi_n(F_n\times A_{n+1})]_n
=\phi(X_n)\subset \phi(X).
 $$
 Therefore, $\widetilde X\subset\phi(X)\subset \widetilde X$.
 The claim follows. 
 
 {\sl Claim 4.}   $\phi$ is a homeomorphism.
 
 \noindent Of course, $\phi$ is continuous and $\phi(X_n)$ is compact and open in $\widetilde X$.
 Since 
 $X=\bigcup_{n=1}^\infty X_n$ and $\widetilde X=\bigcup_{n=1}^\infty\phi(X_n)$,
 the mapping $\phi^{-1}$ is continuous on $\phi(X_n)$ for each $n$.
 Hence $\phi^{-1}$ is continuous on $\widetilde X$.

 {\sl Claim 5.} $\phi (\Cal R(x))=\widetilde{\Cal R}(\phi(x))$ for each $x\in X$.
 
 \noindent
 This claim follows straightforward from the definition of $\phi$.
 
  {\sl Claim 6.} $\mu\circ\phi^{-1}=\widetilde\mu$.
  
  \noindent
  By Claims~4 and 5, $\mu\circ\phi^{-1}$ is an $\widetilde{\Cal R}$-invariant Radon measure, and hence it is a positive scalar multiple of $\widetilde\mu$.
  Since $a_1=1$, the map $\phi_0:F_0\times A_1\to\widetilde F_0$ is onto, and the definition of $\phi$ gives $\phi(X_0)=\widetilde X_0$.
  The normalizations $\mu(X_0)=\widetilde\mu(\widetilde X_0)=1$ therefore force the scalar to be one.

%It follows from Claim 5 and    Radon unique ergodicity of $\widetilde{\Cal R}$ that $c\widetilde\mu=\mu\circ\phi^{-1}$ for some $c>0$.
 % As $\widetilde X_{n-1}\subset\phi([X_n]_n)\subset\widetilde X_{n}$, we obtain that
 % $\mu(X_n)\ge c\widetilde\mu(\widetilde X_{n-1})$ for each $n$.
 % Hence, $\mu(X)\ge c\widetilde\mu(\widetilde X)$.
  %Now (ii) yields that $c=1$.
  
  Thus, the theorem is proved completely.
 \qed
 \enddemo

 \proclaim{Corollary 2.2} Suppose that $T=(T_g)_{g\in G}$ and $\widetilde T=(\widetilde T_g)_{g\in G}$ are the topological $(C,F)$-actions 
 on $X$ and $\widetilde X$ associated with 
 $\Cal T:=(C_n,F_{n-1})_{n=1}^\infty$ and $\widetilde{\Cal T}:=(\widetilde C_n,\widetilde F_{n-1})_{n=1}^\infty$
respectively. 
Then under the conditions of Theorem~2.1, $T$ and $\widetilde T$ are topologically orbit equivalent.
The corresponding topological orbit equivalence is constructed in an explicit algorithmic way from $\Cal T$ and $\widetilde{\Cal T}$.
 \endproclaim

 \proclaim{Corollary 2.3} Suppose that $T=(T_g)_{g\in G}$ and $\widetilde T=(\widetilde T_g)_{g\in G}$ are the measure-preserving $(C,F)$-actions 
 on $(X,\mu)$ and $(\widetilde X,\widetilde\mu)$ associated with 
 $\Cal T:=(C_n,F_{n-1})_{n=1}^\infty$ and $\widetilde{\Cal T}:=(\widetilde C_n,\widetilde F_{n-1})_{n=1}^\infty$
respectively. 
Then under the conditions of Theorem~2.1, $T$ and $\widetilde T$ are 
(measure-theoretically) orbit equivalent.
The corresponding measurable orbit equivalence is constructed 
in an explicit algorithmic way from $\Cal T$ and $\widetilde{\Cal T}$.
 \endproclaim
 
 Our next aim is to prove Theorem B.

\proclaim{Lemma 2.4} Let $(p_n)_{n=0}^\infty$ and $(q_n)_{n=0}^\infty$ be sequences of natural numbers such that $p_0=q_0=1$, $p_n\mid p_{n+1}$ and $q_n\mid q_{n+1}$ for every $n$.
Let $(u_n)_{n=0}^\infty$ and $(v_n)_{n=0}^\infty$ be nondecreasing sequences of positive real numbers with a common limit $L\in(0,+\infty]$, and suppose that $u_n<L$ and $v_n<L$ for every $n$.
If
$$
\bigcup_{n\ge0}\frac1{p_n}\Bbb Z=\bigcup_{n\ge0}\frac1{q_n}\Bbb Z,
\tag{2-3}
$$
then there are strictly increasing sequences of integers
$$
0=k_0<k_1<\cdots\quad\text{and}\quad 0=l_0<l_1<\cdots
$$
such that, for every $n\ge1$,
$$
q_{l_{n-1}}\mid p_{k_n}\mid q_{l_n}
\quad\text{and}\quad
v_{l_{n-1}}\le u_{k_n}\le v_{l_n}.
\tag{2-4}
$$
The two sequences can be selected algorithmically.
\endproclaim

\demo{Proof}
Suppose that $k_{n-1}$ and $l_{n-1}$ have already been chosen.
By~\thetag{2-3}, $1/q_{l_{n-1}}$ belongs to the group on the left-hand side, and hence $q_{l_{n-1}}$ divides $p_j$ for some $j$.
Divisibility then holds for all larger indices.
Since $u_j\nearrow L$ and $v_{l_{n-1}}<L$, we can choose $k_n>k_{n-1}$ so that
$$
q_{l_{n-1}}\mid p_{k_n}
\quad\text{and}\quad
v_{l_{n-1}}\le u_{k_n}.
$$
Applying~\thetag{2-3} in the opposite direction and using $u_{k_n}<L$, we can then choose $l_n>l_{n-1}$ so that
$$
p_{k_n}\mid q_{l_n}
\quad\text{and}\quad
u_{k_n}\le v_{l_n}.
$$
Thus~\thetag{2-4} holds.  Choosing at each step the first pair of indices satisfying these finite divisibility and inequality tests gives the asserted algorithm.
\qed
\enddemo

\demo{Proof of Theorem B}
$(\Longleftarrow)$ We proceed in 3 steps.

Step {\bf I.} Suppose first that $c=1$.
Write
$$
p_n:=\#C_1\cdots\#C_n,\qquad
q_n:=\#\widetilde C_1\cdots\#\widetilde C_n,
$$
with $p_0=q_0=1$, and put
$$
u_n:=\mu(X_n)=\frac{\#F_n}{p_n},
\qquad
v_n:=\widetilde\mu(\widetilde X_n)=\frac{\#\widetilde F_n}{q_n}.
$$
The two sequences of measures increase to $\mu(X)$ and $\widetilde\mu(\widetilde X)$ respectively.
Because $X$ and $\widetilde X$ are non-compact, every $X_n$ and $\widetilde X_n$ is a proper compact open subset of the corresponding space.
Since the Haar measures have full support,
$$
u_n<\mu(X)\quad\text{and}\quad v_n<\widetilde\mu(\widetilde X)
$$
for every $n$.
The assumption of Theorem~B and Lemma~2.4 now give strictly increasing sequences $(k_n)_{n=0}^\infty$ and $(l_n)_{n=0}^\infty$ satisfying~\thetag{2-4}.

Let $\Cal T'$ be the $\boldsymbol k$-telescoping of $\Cal T$ and let $\widetilde{\Cal T}'$ be the $\boldsymbol l$-telescoping of $\widetilde{\Cal T}$.
Thus
$$
C_n'=C_{k_{n-1}+1}\cdots C_{k_n},\qquad F_n'=F_{k_n},
$$
and similarly for the tilded parameters.
For $n\ge1$, set
$$
a_n:=\frac{q_{l_{n-1}}}{p_{k_{n-1}}}
\quad\text{and}\quad
b_n:=\frac{p_{k_n}}{q_{l_{n-1}}}.
$$
These are natural numbers by~\thetag{2-4}, and $a_1=1$.
Moreover,
$$
\#C_n'=\frac{p_{k_n}}{p_{k_{n-1}}}=a_nb_n\qquad\text{and}\qquad
\#\widetilde C_n'=\frac{q_{l_n}}{q_{l_{n-1}}}=b_na_{n+1}.
$$
The  inequalities in~\thetag{2-4} become
$$
\widetilde\mu(\widetilde X_{l_{n-1}})
\le \mu(X_{k_n})
\le \widetilde\mu(\widetilde X_{l_n}),
$$
which are precisely condition~(ii) of Theorem~2.1 for the two telescoped sequences.
Therefore, Corollary~2.2 provides an explicit homeomorphism between the telescoped $(C,F)$-spaces which maps every orbit onto an orbit.
We note that the telescoping mappings from \S1.3 are explicit topological conjugacies.
Composing them with the homeomorphism furnished by Corollary~2.2 (or Theorem~2.1) gives the required map $\phi:X\to\widetilde X$.
The measure identity follows from Theorem~2.1.
All choices can be made by taking the first admissible level and, after fixing enumerations of the finite parameter sets, the first admissible finite injection at every stage; hence the construction is algorithmic in the parameters.

Step {\bf  II.} Let $c$ be a proper  divisor of 
 $\#C_1\cdots\#C_N$ for some $N\in\Bbb N$.
 Select a subset $A$ in $C_1\cdots C_N$ such that $\# A=\#C_1\cdots\#C_N/c$.
We now put
$$
\gather
\text{$ F_0':=F_0$, $ F_1' := F_N$ and $F_n' := F_{n+N-1}$ if  $n\geq 2$ and}\\ 
\text{$ C_1' := A$ and
$ C_n' := C_{n+N-1}$ if $n \geq 2$.}
\endgather
$$
It is routine to verify that 
 $ {\Cal T}':=( C_n', F_{n-1}')_{n=1}^\infty$ 
 is a well defined $(C,F)$-sequence
 that determines a topological $(C,F)$-action $ T'$ on a locally compact Cantor space $X'$.
 Of course, $X'_n= X_{N+n-1}$ for all $n\ge 1$.
 Hence,
there is a canonical identification $\zeta:X\to X'$  of $X$ with $X'$
that intertwines $T$ with $T'$ and
 $\mu\circ\zeta^{-1}=\frac1c\cdot \mu'$, where $\mu'$ is the Haar measure on $ X'$.
Since $\#  C_1' = (\#C_1 \#C_2\ldots \#C_N) / c$, it follows that  $\Cal D({\Cal T}') = c \cdot  \Cal D(\Cal T)$ and $\mu'\circ\zeta = c \cdot \mu$.
Hence, $\mu'(X') = c \cdot \mu(X)$.

Step {\bf III.}  Consider the general case, i.e. let $c$ be arbitrary.
Since $\Cal D(\Cal T)\cup \Cal D(\widetilde{\Cal T})\subset\Bbb Q$, we have that $c\in\Bbb Q$.
Hence, there exist $d,e\in\Bbb N$ such that 
$$
c=d/e,\quad  (d,e)=1,\quad
e \cdot \Cal D(\Cal T) = d \cdot  \Cal D(\widetilde{\Cal T})\quad\text{and}\quad e\cdot\mu(X)=d\cdot\widetilde \mu(\widetilde X).
$$
It follows that the group $e \cdot \Cal D(\Cal T)$ contains $e$ and $d$.
Since $e$ and $d$ are coprime, $1\in e \cdot \Cal D(\Cal T)$.
This implies, in turn, that there is $N$ such that $e$ is a proper divisor  of $\#C_1\cdots\#C_N$.
By Step II, there is a canonical  identification $\zeta:X\to X'$ with $(\zeta T_g\zeta^{-1})_{g\in G}=T'$, $\mu'(X')=e\cdot\mu(X)$ and 
 $\Cal D({\Cal T}') = e \cdot  \Cal D(\Cal T)$.
 In a similar way, 
 there is a canonical  identification $\widetilde \zeta:\widetilde X\to \widetilde X'$ with $(\widetilde \zeta \widetilde T_g\widetilde \zeta^{-1})_{g\in G}=\widetilde T'$, $\widetilde \mu'(\widetilde X')=d\cdot\widetilde \mu(\widetilde X)$ and 
 $\Cal D(\widetilde {\Cal T}') = d \cdot  \Cal D(\widetilde {\Cal T})$.
It remains to apply Step~I to $T'$ and $\widetilde T'$.

$(\Longrightarrow)$ We first show that
$$
\{\mu(K)\mid K\subset X\text{ is compact and open}\}
=
\Cal D(\Cal T)\cap[0,\mu(X)).
\tag{2-5}
$$
Every compact open subset of $X$ is, at some level $n$, a finite union of $n$-cylinders, and each $n$-cylinder has measure $1/p_n(\Cal T)$.
This proves 
%The inclusion 
``$\subset$'' in~\thetag{2-5}.
% is obvious.
Conversely, take $d=m/p_n(\Cal T)\in\Cal D(\Cal T)$ with $0\le d<\mu(X)$.
For some $N\ge n$ we have $d\le\mu(X_N)=\#F_N/p_N(\Cal T)$.
Since $p_n(\Cal T)$ divides $p_N(\Cal T)$, the number $d p_N(\Cal T)$ is an integer not exceeding $\#F_N$.
A union of that many $N$-cylinders has measure $d$.
The endpoint $\mu(X)$ is not attained because a compact open subset is proper in the non-compact space $X$, while $\mu$ has full support.
This proves~\thetag{2-5}.

Now let $\phi$ be an orbit equivalence of $T$ with $\widetilde T$  (see the statement of Theorem~B).
The measure $\mu\circ\phi^{-1}$ is invariant under the %full group of the 
$\widetilde T$-orbit equivalence relation and hence, by the Radon unique ergodicity, equals $c\widetilde\mu$ for some $c>0$.
Thus $\phi$ carries the set on the left-hand side of~\thetag{2-5} onto $c$ times the analogous set for $\widetilde T$.
Taking the additive groups generated by these compact-open values gives
$$
\Cal D(\Cal T)=c\cdot \Cal D(\widetilde{\Cal T})\qquad\text{and}\qquad \mu\circ\phi^{-1}=c\cdot\widetilde\mu,
$$
as desired.
\qed
\enddemo

\comment
The converse to Theorem B is also true.

\proclaim{Proposition 2.5} Let $T$ be a continuous $(C,F)$-action on a non-compact locally compact Cantor space $X$, associated with $\Cal T=(C_n,F_{n-1})_{n=1}^\infty$, and let $\mu$ be the Haar measure on $X$.
Then

Consequently, if $T$ and $\widetilde T$ are topologically orbit equivalent via a homeomorphism $\phi:X\to\widetilde X$, then for some $c>0$,
$$
\Cal D(\Cal T)=c\cdot \Cal D(\widetilde{\Cal T})\qquad\text{and}\qquad \mu\circ\phi^{-1}=c\cdot\widetilde\mu.
\tag{2-6}
$$
%If one of the two %invariant Radon measures
%Haar measures (on $X$ and $\widetilde X$)
 %is finite, then so is the other and
%$$
%\mu(X)=c\cdot \widetilde\mu(\widetilde X).
%\tag{2-7}
%$$
\endproclaim

\demo{Proof}

Now suppose that $\phi$ is an orbit equivalence.
The measure $\mu\circ\phi^{-1}$ is invariant under the %full group of the 
$\widetilde T$-orbit equivalence relation and hence, by the Radon unique ergodicity, equals $c\widetilde\mu$ for some $c>0$.
Thus $\phi$ carries the set on the left-hand side of~\thetag{2-5} onto $c$ times the analogous set for $\widetilde T$.
Taking the additive groups generated by these compact-open values gives~\thetag{2-6}.
%Equation~\thetag{2-7} follows by evaluating the two measures on the whole space whenever the total mass is finite.
\qed
\enddemo

\endcomment

We now give an example of two non-topologically orbit equivalent continuous $(C,F)$-transformations on finite measure spaces which are non-compact but locally compact.
Moreover, the spaces  are of the  same finite Haar measure.

\example{Example 2.5}
We use additive notation for $G=\Bbb Z$.
For $r\in\{2,3\}$, define a $(C,F)$-sequence inductively as follows.
Let $F_0^{(r)}=\{0\}$.
If
$$
F_n^{(r)}=[s_n^{(r)},t_n^{(r)}]\cap\Bbb Z
\quad\text{and}\quad h_n^{(r)}:=\#F_n^{(r)},
$$
we put
$$
\align
C_{n+1}^{(r)}&:=\{0,h_n^{(r)},2h_n^{(r)},\dots,(r-1)h_n^{(r)}\},\\
F_{n+1}^{(3)}&:=
[s_n^{(3)}-1,\ t_n^{(3)}+2h_n^{(3)}+1]\cap\Bbb Z,\\
F_{n+1}^{(2)}&:=
[s_n^{(2)}-1,\ t_n^{(2)}+h_n^{(2)}]\cap\Bbb Z\quad \text{if $n$ is odd and}\\
F_{n+1}^{(2)}&:=
[s_n^{(2)},\ t_n^{(2)}+h_n^{(2)}+1]\cap\Bbb Z\quad \text{if $n$ is even.}\\
\endalign
$$
The translates $F_n^{(r)}+c$, $c\in C_{n+1}^{(r)}$, are pairwise disjoint and contained in $F_{n+1}^{(r)}$.
%Furthermore, the set
%$$
%F_n^{(r)}+C_{n+1}^{(r)}+\cdots+C_m^{(r)}
%$$
%has distance at least $m-n$ from both endpoints of $F_m^{(r)}$.
It is  routine to verify that
the condition~\thetag{1-4} holds, so the associated $(C,F)$-action of $\Bbb Z$ is continuous on the entire $(C,F)$-space $X^{(r)}$.
Thus, we add  either one or two spacers at every step  of the $(C,F)$-construction, and therefore $X^{(r)}$ is non-compact.
Let $\mu^{(r)}$ denote the Haar measure on $X^{(r)}$.
%Since $h_{n+1}^{(r)}=rh_n^{(r)}+2,$
We have that
$$
\mu^{(3)}(X^{(3)})=1+\sum_{n=1}^\infty\frac{2}{3^n}=2\qquad\text{and}\qquad
\mu^{(2)}(X^{(2)})=1+\sum_{n=1}^\infty\frac{1}{2^n}=2.
$$
Thus, $\mu^{(3)}(X^{(3)})=\mu^{(2)}(X^{(2)})$.
On the other hand,
$$
\Cal D(\Cal T^{(2)})=\Bbb Z[1/2]
\quad\text{and}\quad
\Cal D(\Cal T^{(3)})=\Bbb Z[1/3].
$$
%No positive scalar multiple of $\Bbb Z[1/3]$ equals $\Bbb Z[1/2]$: the former is divisible by $3$, whereas the latter is not.
Hence, Proposition~2.5 shows that these two continuous rank-one actions are not topologically orbit equivalent.
In this connection we also note that no positive scalar multiple of $\Cal D(\Cal T^{(3)})$ equals $\Cal D(\Cal T^{(2)})$: the former is divisible by $3$, whereas the latter is not.
%Thus the unrestricted formulation of Theorem~B is false.
\endexample

\head 3. Measurable orbit equivalence and Proof of Theorem A
\endhead

We start this section with an auxiliary lemma.

\proclaim{Lemma 3.1} Let $(x_n)_{n=1}^\infty$ and $(y_n)_{n=1}^\infty$ be a sequence of positive reals and a sequence of natural numbers respectively such that  $y_n\to\infty$ and $\lim_{n\to\infty}\frac{x_n}{y_n}=\alpha\in (0,+\infty)$. 
Then for each $\epsilon>0$, $L>\alpha$ and $d\in\Bbb N$, there exists $N$ such that  
for each $n>N$, the finite set
$$
\bigg\{\frac {x_n}{i}\,\bigg|\,  i\in\Bbb N, i<y_n\text{ and } \,d|i\bigg\}\cap [\alpha,L]
$$
is an $\epsilon$-net in $[\alpha, L]$.
\comment

$$
\Big|\frac{x_n}{y_n-di}-\frac{x_n}{y_n-d(i-1)}\Big|<\epsilon
$$
 for each  $n>N$ and every $i\in\Bbb N $
with  $\frac{x_n}{y_n-di}< L$. 
\endcomment
\endproclaim

We omit the proof because it is  routine.

\proclaim{Lemma 3.2} Let $(f_n)_{n=1}^\infty$, $(c_n)_{n=1}^\infty$,
$(\widetilde f_n)_{n=1}^\infty$, $(\widetilde c_n)_{n=1}^\infty$ be  sequences
of natural numbers such that  $c_n\ge 2$, $\widetilde c_n\ge 2$ for each $n$ and  there exist  
$\lim_{n\to\infty}\frac{f_n}{c_1\cdots c_n}=\alpha\in (0,\infty]$ and 
$\lim_{n\to\infty}\frac{\widetilde f_n}{\widetilde c_1\cdots \widetilde c_n}
=\widetilde \alpha\in (0,\infty]$.
If $\max(\alpha,\widetilde\alpha)<\infty$ then there are two increasing sequences 
$(k_n)_{n=0}^\infty$ and 
$(l_n)_{n=0}^\infty$
and two sequences $(a_n)_{n=1}^\infty$ and $(b_n)_{n=1}^\infty$ such that $k_n,l_n,a_n,b_n\in\Bbb N$, $k_0=l_0=0$,
$a_1=1$,  
$$
\gather
2\le a_nb_n<c_{k_{n-1}+1}\cdots c_{k_n}, \quad 2\le b_na_{n+1}<\widetilde c_{l_{n-1}+1}\cdots \widetilde c_{l_n}\quad\text{and }
\tag3-1\\
\frac{f_{k_n}}{a_1b_1\cdots a_nb_n}<\frac{\widetilde f_{l_n}}{b_1a_2\cdots b_na_{n+1}}<
\frac{f_{k_{n+1}}}{a_1b_1\cdots a_{n+1}b_{n+1}}<2\max(\alpha,\widetilde\alpha).
\tag3-2
\endgather
$$
for each $n\in\Bbb N$.
If $\alpha=\widetilde\alpha=\infty$ then the same holds but
with 
$$
\frac{a_nb_n}{c_{k_{n-1}+1}\cdots c_{k_n}}>1-\frac 1{n^2}\quad \text{ and }\quad
\frac{b_na_{n+1}}{\widetilde c_{l_{n-1}+1}\cdots \widetilde c_{l_n}}>1-\frac1{n^2}
$$
in place of the rightmost inequality in  \thetag{3-2}.
\endproclaim

\demo{Proof}
Suppose first that $\beta:=\max(\alpha,\widetilde\alpha)<\infty$.
Let $a_1:=1$.
%Fix a sequence of positive reals $(\epsilon_n)_{n=1}^\infty$ such that $\sum_{n=1}^\infty\epsilon_n<0.1$.
It follows from  Lemma~3.1 that there exist  $k_1>1$ and $b_1\in\{2,3,\dots\}$ such that 
$$
a_1b_1< c_1\cdots c_{k_1},\quad         \frac{f_{k_1}}{a_1b_1}\in \Big(\frac32\beta, 2\beta\Big)
\quad \text{and}\quad   \frac{f_{k_1}}{a_1b_1}\cdot\frac{c_1\cdots c_{k_1}}{f_{k_1}}\cdot\alpha <2\beta.%\quad   \inf_{n>k_1}\frac{a_1b_1}{f_{k_1}}\cdot \frac{f_n}{c_1\cdots c_n}\in (1-\epsilon_1,1+\epsilon_1).
\tag3-3
$$
For instance, choose $k_1$  large so that 
\roster
\item"---" $\frac{c_1\cdots c_{k_1}}{f_{k_1}}\cdot\alpha <1.1$
and 
\item"---" there is $b_1$  satisfying the following inequalities: $\frac32\beta\le \frac{f_{k_1}}{a_1b_1}\le\frac53\beta$.
\endroster
%We have that
%$$
%\frac{\widetilde f_n}{\widetilde c_1\cdots\widetilde c_{n}}=
%\frac{\widetilde f_n}{\widetilde c_1\cdots\widetilde c_{n}}\frac{\widetilde f_n}{\widetilde c_1\cdots\widetilde c_{n}}
%$$
In a similar way, applying Lemma~3.1 again,  we find $l_1>1$ and $a_2\in\Bbb N$ such that 
$$
b_1a_2< \widetilde c_1\cdots \widetilde c_{l_1},  \ 
\frac{\widetilde f_{l_1}}{b_1a_2}\in \bigg(\frac{f_{k_1}}{a_1b_1}, 2\beta\bigg)
\quad\text{and}\quad
\frac{\widetilde f_{l_1}}{b_1a_2}\cdot\frac{\widetilde c_1\cdots \widetilde c_{l_1}}{\widetilde f_{l_1}}\cdot\widetilde\alpha <2\beta.
\tag3-4
$$
%for each $j=0,1,\dots,b_1-1$.
%Choosing   $j$ so that $b_1|(\widetilde c_1\cdots \widetilde c_{l_1}-j)$, we obtain that there is $a_2$ such that
%$b_1a_2<\widetilde c_1\cdots\widetilde c_{l_1}$ and 
%$$
%\frac{\widetilde f_{l_1}}{b_1a_2}\in \bigg(\frac{f_{k_1}}{a_1b_1}, 2\beta\bigg).
%$$
It follows from the rightmost  inequality in \thetag{3-3} and the assumption of the lemma that
$$
%\align
\lim_{n\to\infty}\frac{f_n}{a_1b_1c_{k_1+1}\cdots c_n}= \frac{c_1\cdots c_{k_1}}{a_1b_1} \cdot\alpha<2\beta.
%\\
%&\in \Big(\frac32\beta, 2\beta\Big)\cdot (1-\epsilon_1,1+\epsilon_1)\\
%&=\Big(\frac32\beta(1-\epsilon_1), 2\beta(1+\epsilon_1)\Big)\\
%&\ni 
%\frac{\widetilde f_{l_1}}{b_1a_2}.
%\endalign
$$
Hence, by Lemma~3.1,  there exist $k_2>k_1$ and $b_2\in\{2,3,\dots\}$ such that 
$$
a_2b_2< c_{k_1+1}\cdots c_{k_2}, \   
\frac{ f_{k_2}}{a_1b_1a_2b_2}\in \bigg(\frac{\widetilde f_{l_1}}{b_1a_2},2\beta\bigg)
\quad \text{and}\quad   %\inf_{n>k_2}
\frac{f_{k_2}}{a_1b_1a_2b_2}\cdot \frac{c_1\cdots c_{k_2}}{f_{k_2}}\cdot\alpha<2\beta.
$$
%for each $j=0,1,\dots, b_1$
In a similar way, we deduce from the rightmost inequality in \thetag{3-4} that
$$
\lim_{n\to\infty}\frac{\widetilde f_n}{b_1a_2\widetilde c_{l_1+1}\cdots\widetilde  c_n}=
 \frac{\widetilde c_1\cdots\widetilde c_{l_1}}{b_1a_2} \cdot\widetilde \alpha<2\beta.
$$
Continuing this by induction in $n$, we construct increasing sequences
$(k_n)_{n=0}^\infty$ and 
$(l_n)_{n=0}^\infty$
and
sequences  $(a_n)_{n=1}^\infty$ and $(b_n)_{n=1}^\infty$
such that $k_n,l_n,a_n,b_n$ are positive integers if $n>0$, $k_0=l_0=0$, $a_1=1$ and \thetag{3-1} and \thetag{3-2} hold for each $n$.

If $\alpha=\widetilde\alpha=\infty$ then the argument is easier as one needs not to apply Lemma~3.1.
We leave details to the reader.
\qed
\enddemo

We note that  \thetag{3-1}, \thetag{3-2} and the hypotheses of Lemma~3.2 yield that
$$
0<\prod_{n=1}^\infty\frac{a_nb_n}{c_{k_{n-1}+1}\cdots c_{k_n}}<\infty\quad\text{and}\quad
0<\prod_{n=1}^\infty\frac{b_na_{n+1}}{\widetilde c_{l_{n-1}+1}\cdots \widetilde c_{l_n}}<\infty.
\tag3-5
$$
This is obvious if  $\alpha=\widetilde\alpha=\infty$.
If $\beta<\infty$ then for each $j>1$,
$$
\prod_{n=1}^j\frac{a_nb_n}{c_{k_{n-1}+1}\cdots c_{k_n}}=\frac{f_{k_j}/(c_1\cdots c_{k_j})}{f_{k_j}/(a_1b_1\cdots a_jb_j)}.
$$
The numerator tends to $\alpha$, while the denominator
increases  as $j\to\infty$ and  is bounded  above by $2\beta$.
Hence the ratio has a finite positive limit.
A similar argument applies to $\prod_{n=1}^\infty\frac{b_na_{n+1}}{\widetilde c_{l_{n-1}+1}\cdots \widetilde c_{l_n}}$.
\comment

\proclaim{Lemma 2.1}
Let $(b_n)_{n=1}^\infty$ be a sequence of natural numbers such that $\sum_{n=1}^\infty\frac 1{b_n}<\infty$. 
Let $\beta:=\prod_{n=1}^\infty\big(1-\frac 1{b_n}\big)$.
Then for each $\alpha\in (0,\beta)$, there exists a sequence $(a_n)_{n=1}^\infty$ such that
$0<a_n<b_n$ for each $n$ and $\alpha=\prod_{n=1}^\infty\big(1-\frac {a_n}{b_n}\big)$.
\endproclaim

\proclaim{Corollary 2.2}
Let $(b_n)_{n=1}^\infty$ and $(\widetilde b_n)_{n=1}^\infty$ be two sequences of natural numbers such that $\sum_{n=1}^\infty\frac 1{b_n}<\infty$ and
$\sum_{n=1}^\infty\frac 1{\widetilde b_n}<\infty$. 
Let 
$$
\beta:=\prod_{n=1}^\infty\big(1-\frac 1{b_n}\big)\text{ and }\widetilde\beta:=\prod_{n=1}^\infty\big(1-\frac 1{\widetilde b_n}\big).
$$
Then for each  $\gamma\in (0,\min(\beta,\widetilde\beta))$,
there exist  two sequences $(a_n)_{n=1}^\infty$  and $(\widetilde a_n)_{n=1}^\infty$ such that
$0<a_n<b_n$, $0<\widetilde a_n<\widetilde b_n$ and
$$
\alpha=\prod_{n=1}^\infty\big(1-\frac {a_n}{b_n}\big)=\prod_{n=1}^\infty\big(1-\frac {\widetilde a_n}{\widetilde b_n}\big).
$$
\endproclaim

\proclaim{Theorem 3.3} Given two  $(C,F)$-actions which  either both preserve finite measure or both preserve  infinite measure, there is an explicit algorithmic construction
 of a measurable orbit equivalence between these  actions.
\endproclaim

\endcomment

\demo{Proof of Theorem A}
%Let $(X,\mu,T)$ and $(\widetilde X,\widetilde\mu,\widetilde T)$ be
% two measure preserving $(C,F)$-actions 
Let $T$ and $\widetilde T$ be associated with
 $(C,F)$-sequences $\Cal T:=(C_n,F_{n-1})_{n=1}^\infty$ and $\widetilde{\Cal T}:=(\widetilde C_n,\widetilde F_{n-1})_{n=1}^\infty$
respectively.
%Suppose that $\mu(X)<\infty$ and $\widetilde\mu(\widetilde X)<\infty$. 
Let $f_n:=\# F_n$, $c_n:=\#C_n$, $\widetilde f_n:=\# \widetilde F_n$ and  $\widetilde c_n:=\#\widetilde C_n$.
Since 
$$
\mu(X)=\lim_{n\to\infty}\frac{\#F_n}{\#C_1\cdots\# C_n}
\quad
\text{and}\quad
\widetilde\mu(\widetilde X)=\lim_{n\to\infty}\frac{\#\widetilde F_n}{\#\widetilde C_1\cdots\# \widetilde C_n},
$$  
we may apply Lemma~3.2 and find increasing sequences
$\boldsymbol k:=(k_n)_{n=0}^\infty$ and 
$\boldsymbol l:=(l_n)_{n=0}^\infty$
and two sequences $(a_n)_{n=1}^\infty$ and $(b_n)_{n=1}^\infty$ such that $k_0=l_0=0$, $a_1=1$ and~\thetag{3-1} and~\thetag{3-2} (without the rightmost inequality) hold for each $n\in\Bbb N$.
In view of \thetag{3-1}, we can choose subsets $C_n^\circ\subset C_{k_{n-1}+1}\cdots C_{k_n}$ and 
$\widetilde C_n^\circ\subset \widetilde C_{l_{n-1}+1}\cdots \widetilde C_{l_n}$ such that $\# C_n^\circ=a_nb_n$ and $\#\widetilde C_n^\circ=b_na_{n+1}$ for each $n$.
According to \thetag{3-5},
$$
0<\prod_{n=1}^\infty\frac{\# C_n^\circ}{\# C_{k_{n-1}+1}\cdots \# C_{k_n}}<\infty\quad\text{and}\quad
0<\prod_{n=1}^\infty\frac{\#\widetilde C_n^\circ}{\#\widetilde C_{l_{n-1}+1}\cdots \#\widetilde C_{l_n}}<\infty.\tag3-6
$$
Let $\boldsymbol C:=(C_n^\circ)_{n=1}^\infty$ and 
$\widetilde{\boldsymbol C}:=(\widetilde C_n^\circ)_{n=1}^\infty$.
Denote by $\Cal T'$ the $\boldsymbol C$-reduction of the $\boldsymbol k$-telescoping of
$\Cal T$
and denote by $\widetilde{\Cal T}'$
the $\widetilde{\boldsymbol C}$-reduction  of the $\boldsymbol l$-telescoping of $\widetilde{\Cal T}$.
These reductions are well defined due to \thetag{3-6}.
Then  $\Cal T'$ and $\widetilde{\Cal T}'$ satisfy the conditions of Theorem~2.1.
Indeed, (i) follows from \thetag{3-1} and~(ii) follows from \thetag{3-2}.
Hence, by Theorem~2.1, the $(C,F)$-actions $T'$ and $\widetilde T'$ associated with  $\Cal T'$ and $\widetilde{\Cal T}'$ 
are measurably orbit equivalent. 
As $(X,\mu,T)$ and $(\widetilde X,\widetilde\mu,\widetilde T)$ are isomorphic (via measure scaling mappings)
to $T'$ and $\widetilde T'$ respectively, it follows that
$T$ and $\widetilde T$ are also
measurably orbit equivalent.
More precisely, 
let $\iota_{\boldsymbol k}$ and $\iota_{\boldsymbol l}$ denote the $\boldsymbol k$-telescoping mapping and the $\boldsymbol l$-telescoping mapping respectively.
Let $\rho_{\boldsymbol C}$ and $\rho_{\widetilde{\boldsymbol C}}$
stand for the $\boldsymbol C$-reduction mapping and $\widetilde{\boldsymbol C}$-reduction mapping respectively.
We remind that $\iota_{\boldsymbol k}$ and $\iota_{\boldsymbol l}$ are homeomorphisms
and that 
$\rho_{\boldsymbol C}$ and $\rho_{\widetilde{\boldsymbol C}}$ are Borel and 
defined on conull invariant dense $F_\sigma$-subsets of the corresponding
$(C,F)$-spaces.
Denote these subsets by $Y$ and $\widetilde Y$ respectively.
Let $\psi$ be the homeomorphism of $\rho_{\boldsymbol C}(Y)$ onto 
$\rho_{\widetilde{\boldsymbol C}}(\widetilde Y)$
 supplied by Theorem~2.1.
 We put
$$
\align
X^\circ&:=X^G\cap \iota_{\boldsymbol k}^{-1}\Big(Y\cap \rho^{-1}_{\boldsymbol C}\circ\psi^{-1}\circ\rho_{\widetilde{\boldsymbol C}}\big(\widetilde Y\cap\iota_{\boldsymbol l}(\widetilde X^G)\big)\Big), \\
\widetilde X^\circ&:= \widetilde X^G\cap 
\iota_{\boldsymbol l}^{-1}\Big(\widetilde Y\cap \rho^{-1}_{\widetilde{\boldsymbol C}}\circ\psi\circ\rho_{\boldsymbol C}\big( Y\cap\iota_{\boldsymbol k}(X^G)\big)\Big)\quad\text{and}\\
\phi &:=\iota_{\boldsymbol l}^{-1}\circ\rho_{\widetilde{\boldsymbol C}}^{-1}\circ\psi\circ\rho_{\boldsymbol C}\circ\iota_{\boldsymbol k}.
\endalign
$$
Then $X^\circ$ and $\widetilde X^\circ$ are conull invariant %dense $F_\sigma$-
subsets of $X^G$ and 
$\widetilde X^G$ respectively
and
$\phi:X^\circ\to\widetilde X^\circ$ is a well defined measurable orbit equivalence of $T$ and $\widetilde T$.
Hence,  there is  $c>0$ such that $\widetilde\mu\circ\phi=c\mu$.
If $\mu(X)+\widetilde\mu(\widetilde X)<\infty$ then
$c=\frac{\widetilde\mu(\widetilde X)}{\mu(X)}$.
If $\mu(X)=\widetilde\mu(\widetilde X)=\infty$ then 
$$
c=\prod_{n=1}^\infty\left(\frac{\# C_{k_{n-1}+1}\cdots \#C_{k_n}}{\# C_n^\circ}\cdot
\frac{\# \widetilde C_n^\circ}{\# \widetilde C_{l_{n-1}+1}\cdots \#\widetilde C_{l_n}}\right).
$$
\qed
 %The case where $\mu(X)=\widetilde\mu(\widetilde X)=\infty$ is considered in a similar way.
\enddemo

\comment

\example{Example 3.4 {\rm (cf. \cite{HaKe})}} Let $T$ and $\widetilde T$ be the 2-odometer and the 3-odometer respectively. 
Let $F_n=\{0,\dots, 2^n-1\}$, $C_{n+1}=\{0,2^{n}\}$,
$\widetilde F_n=\{0,\dots, 3^n-1\}$, $\widetilde C_{n+1}=\{0,3^{n}, 2\cdot 3^n\}$.
Then $T$ and $\widetilde T$ are the $(C,F)$-transformations ($\Bbb Z$-actions) associated with 
 $\Cal T:=(C_n,F_{n-1})_{n=1}^\infty$ and $\widetilde{\Cal T}:=(\widetilde C_n,\widetilde F_{n-1})_{n=1}^\infty$ respectively.
We note that
$\#F_n=2^n$, $\#C_{n+1}=2$, $\#\widetilde F_n=3^n$ and $\#\widetilde C_{n+1}=3$ for all $n\ge 0$.
We will construct inductively $k_n,l_n,a_n,b_n$ for $n=0,1,\dots$
as in the proof of  Theorem~3.3.
On the first step we set $k_0:=l_0:=0$, $a_1:=1$. 
Suppose that we have defined $k_{n-1},l_{n-1},b_{n-1},a_n$ for some $n\ge 1$ and
$$
1\le\frac{\# \widetilde F_{l_{n-1}}}{a_n\# F_{k_{n-1}}}<\frac{2\prod_{j=1}^{n-1}(a_jb_j)}{\# F_{k_{n-1}}}.
$$
%Suppose that $a_n$ is odd.
Denote by $\phi(.)$  the Euler function.
Then, by the Fermat-Euler theorem, 
 for each $r\in\Bbb N$,
 there exists a number $b(r)$ such that
$$
2^{\phi(3ra_n)}=a_n\cdot 3rb(r)+1=a_n\cdot (3rb(r)-j)+ ja_n+1,
$$
for each $j\in\Bbb N$.
We note that 
$$
1<\frac{2^{\phi(3ra_n)}}{a_n\cdot 3rb(r)}< \frac{2^{\phi(3ra_n)}}{a_n\cdot (3rb(r)-1)}
<\frac{2^{\phi(3ra_n)}}{a_n\cdot (3rb(r)-2)}<\cdots.
$$
Start with $r:=1$ and find
 the smallest integer $j(r)\ge 0$ such that
 $$
 \frac{\# \widetilde F_{l_{n-1}}}{a_n\# F_{k_{n-1}}}\le\frac{2^{\phi(3ra_n)}}{a_n\cdot (3rb(r)-j)}.
 \tag3-7
 $$
If
$$
\frac{2^{\phi(3ra_n)}}{a_n\cdot (3rb(r)-j)}\ge \frac{2\prod_{j=1}^{n-1}(a_jb_j)}{\# F_{k_{n-1}}}
$$
then pass to $r:=2$ and repeat the procedure.
If the right-hand side in \thetag{3-7}  is still greater or equal   
$\frac{2\prod_{i=1}^{n-1}(a_ib_i)}{\# F_{k_{n-1}}}$   then put $r:=3$
and so on. 
After finitely many steps we find $r_n\in\Bbb N$ such that
$$
\frac{\# \widetilde F_{l_{n-1}}}{a_n\# F_{k_{n-1}}}
\le \frac{2^{\phi(3r_na_n)}}{a_n\cdot (3r_nb(r_n)-j(r_n))}
<\frac{2\prod_{i=1}^{n-1}(a_ib_i)}{\# F_{k_{n-1}}}.
\tag3-8
$$
We now let $b_n:=3r_nb(r_n)-j(r_n)$ and 
$k_{n}:=k_{n-1}+\phi(3r_na_n)$.
This completes the first $n$-th semistep.
The second  $n$-th semistep is similar.
It follows from \thetag{3-8}
that
$$
1\le\frac{\# F_{k_{n-1}}\cdot 2^{\phi(3r_na_n)}}{b_n\# \widetilde F_{l_{n-1}}}=
\frac{\# F_{k_{n}}}{b_n\# \widetilde F_{l_{n-1}}}
<\frac{2a_n\prod_{i=1}^{n-1}(a_ib_i)}{\# \widetilde F_{l_{n-1}}}=
\frac{2\prod_{i=1}^{n-1}(b_ia_{i+1})}{\# \widetilde F_{l_{n-1}}}.
$$
Thus,
$$
1\le
\frac{\# F_{k_{n}}}{b_n\# \widetilde F_{l_{n-1}}}
<
\frac{2\prod_{i=1}^{n-1}(b_ia_{i+1})}{\# \widetilde F_{l_{n-1}}}.
$$
Then,
for each $q\in\Bbb N$,
 there exists a number $a(q)$ such that
$$
3^{\phi(2qb_n)}=b_n\cdot 2qa(q)+1=b_n\cdot (2qa(q)-i)+ ib_n+1,
$$
for each $i\in\Bbb N$.
As in the first $n$-th semistep, we can find $q_n$ and $i(q_n)$ such that
$$
\frac{\# F_{k_{n}}}{b_n\# \widetilde F_{l_{n-1}}}
\le 
\frac{3^{\phi(2q_nb_n)}}{b_n\cdot (2q_na(q_n)-i(q_n))}<
\frac{2\prod_{i=1}^{n-1}(b_ia_{i+1})}{\# \widetilde F_{l_{n-1}}}.\tag3-9
$$
Then we set $a_{n+1}:=2q_na(q_n)-i(q_n)$ and $l_{n}:=l_{n-1}+\phi(2q_nb_n)$.
This completes the $n$-th step of the inductive construction.

It follows from \thetag{3-8} that
$$
a_nb_n\le 2^{\phi(3r_na_n)}=2^{k_n-k_{n-1}}=\# C_{k_{n-1}+1}\cdots\# C_{k_n}.
$$
Hence, we can choose a subset $C^\circ_n\subset C_{k_{n-1}+1}+\cdots+ C_{k_n}$ of cardinality $a_nb_n$.
In a similar way, choose a subset $\widetilde C^\circ_n\subset \widetilde C_{l_{n-1}+1}+\cdots+ \widetilde C_{l_n}$ of cardinality $b_na_{n+1}$.
It is possible in view of \thetag{3-9}.
Let $\boldsymbol k:=(k_n)_{n=0}^\infty$,
 $\boldsymbol l:=(l_n)_{n=0}^\infty$,
 $\boldsymbol C^\circ:=(C_n^\circ)_{n=1}^\infty$
and $\widetilde{\boldsymbol C}^\circ:=(\widetilde C_n^\circ)_{n=1}^\infty$.
Let $\Cal L$ be the ${\boldsymbol C}^\circ$-reduction of the $\boldsymbol k$-telescoping of $\Cal T$ and let 
$\widetilde {\Cal L}$ be the $\widetilde{\boldsymbol C}^\circ$-reduction of the $\boldsymbol l$-telescoping of $\widetilde{ \Cal T}$.
Denote by $S$ and $\widetilde S$ the $(C,F)$-transformations associated with 
$\widetilde {\Cal L}$ be the $\widetilde{\boldsymbol C}^\circ$ respectively.
Then $S$ is isomorphic to $T$ and $\widetilde S$ is isomorphic to $\widetilde T$.
Let $(X,\mu)$ and  $(\widetilde X,\widetilde\mu)$ be the spaces of  $S$ and $\widetilde S$
respectively. 
Then we have natural filtrations $X_0\subset X_1\subset \cdots$  and 
$\widetilde X_0\subset \widetilde X_1\subset \cdots$
of $X$ and $\widetilde X$ respectively.
Then \thetag{3-8} is equivalent to the following
$$
\widetilde\mu(\widetilde X_{n-1})\le \mu(X_n)<2.
$$
In a similar way, \thetag{3-9} is equivalent to the following
$$
\mu(X_{n})\le \mu(\widetilde X_n)<2.
$$
Hence, the conditions of  Theorem~2.1 are satisfied.
Then Theorem~2.1 provides an explicit algorithmic construction  of a homeomorphism of $X$ onto $\widetilde X$ which is an orbit equivalence
between  $S$ and $\widetilde S$.

\endexample

\endcomment

\comment

Passing to telescopings of $\Cal T$ and $\widetilde{\Cal T}$ if necessary, we may assume without loss of generality that $\sum_{n=1}^\infty\frac{1}{\# C_n}<\infty$ and  $\sum_{n=1}^\infty\frac{1}{\# \widetilde C_n}<\infty$.
We let $b_n:=\# C_n$ and $\widetilde b_n:=\#\widetilde C_n$.
It follows from  Lemma~2.1 that there exist a large number $\alpha>0$ and subsets $C_n^\circ\subset C_n$
and $\widetilde C_n^\circ\subset \widetilde C_n$ such that
$$
\frac{\mu(X)}\alpha=\prod_{n=1}^\infty\frac{\#C_n^\circ}{\# C_n}\quad\text{ and }\quad
\frac {\widetilde\mu(\widetilde X)}\alpha=\prod_{n=1}^\infty\frac{\#\widetilde C_n^\circ}{\# \widetilde C_n}.\tag2-1
$$
Let $\Cal T'$ and $\widetilde{\Cal T}'$ denote the  $(C_n^\circ)_{n=1}^\infty$-reduction of  $\Cal T$ and the
$(\widetilde C_n^\circ)_{n=1}^\infty$-reduction  of  $\widetilde{ \Cal T}$ respectively.
Let 
Let $(X',\mu',T')$ and $(\widetilde X',\widetilde\mu',\widetilde T')$ be
the  measure-preserving $(C,F)$-dynamical systems associated with
$\Cal T'$ and $\widetilde{\Cal T}'$ respectively.
Since
$$
\frac{\mu(X)}{ \mu'( X')}=\prod_{n=1}^\infty\frac{\#C_n^\circ}{\# C_n}\quad\text{ and }\quad
\frac{\widetilde\mu(\widetilde X)}{\widetilde \mu'(\widetilde X')}=\prod_{n=1}^\infty\frac{\#\widetilde C_n^\circ}{\# \widetilde C_n},
$$
it follows from \thetag{2-1} that $\mu'(X')= \widetilde \mu'( \widetilde X')=\alpha$.
Since the two sequences $(\mu'(X'_n))_{n=1}^\infty$ and $(\widetilde \mu'( \widetilde X_n'))_{n=1}^\infty$
increase and converge to $\alpha$, there exist  two increasing sequences $0=k_0<k_1<\cdots$ and $0=l_0<l_1<\cdots$
such that 
$$
\widetilde \mu'(\widetilde X'_{k_0})\le \mu'( X'_{l_1})\le \widetilde \mu'(\widetilde X'_{k_1})\le
\mu'( X'_{l_2})\le\cdots.
$$
Let $\boldsymbol k:=(k_n)_{n=0}^\infty$ and $\boldsymbol l:=(l_n)_{n=0}^\infty$.
Consider the $\boldsymbol k$-telescoping of $\Cal T'$ and 
 the $\boldsymbol l$-telescoping of $\widetilde{\Cal T}'$.
 Passing, if necessary to  further telescopings and then to reductions we obtain finally 
 two $(C,F)$-sequences that satisfy (i) and (ii) from the statement of Theorem~1.1.
 Hence, we can apply Theorem~1.1.

 \endcomment

\comment

\proclaim{Proposition 3}
Let $b_1<b_2<\cdots$ and $\widetilde b_1<\widetilde b_2<\cdots$ be two increasing sequences of natural numbers such that $\sum_{n=1}^\infty\frac 1{b_n}<\infty$ and
$\sum_{n=1}^\infty\frac 1{\widetilde b_n}<\infty$.
Then there exist four  sequences $(k_n)_{n=0}^\infty$, $(l_n)_{n=0}^\infty$,
$(u_n)_{n=1}^\infty$ and $(v_n)_{n=1}^\infty$ of natural numbers such that
$k_0=0$, $l_0=0$, $u_1=1$, $u_nv_n\le b_{k_{n-1}+1}\cdots b_{k_n}$,
$v_nu_{n+1}\le \widetilde b_{l_{n-1}+1}\cdots \widetilde b_{l_n}$ for each $n$,
$$
\prod_{n=1}^\infty\Big( 1-\frac{u_nv_n}{b_{k_{n-1}+1}\cdots b_{k_n}}\Big)=
\prod_{n=1}^\infty\Big( 1-\frac{v_nu_{n+1}}{\widetilde b_{l_{n-1}+1}\cdots \widetilde b_{l_n}}\Big)
$$
and for each $m\in\Bbb N$,
$$
\prod_{n=1}^m\Big( 1-\frac{u_nv_n}{b_{k_{n-1}+1}\cdots b_{k_n}}\Big)
\le 
\prod_{n=1}^m\Big( 1-\frac{v_nu_{n+1}}{\widetilde b_{l_{n-1}+1}\cdots \widetilde b_{l_n}}\Big)
\le
\prod_{n=1}^{m+1}\Big( 1-\frac{u_nv_n}{b_{k_{n-1}+1}\cdots b_{k_n}}\Big).
$$
\endproclaim

\endcomment

\enddocument